\documentclass[12pt]{article}

\catcode `\@=11
\def\address{\m@th\@ifnextchar[\@address{\@address[]}}

\def\@address[#1]#2{
\expandafter\def\expandafter\@addressname\expandafter
{\@addressname{
  \adr{#1}\ \parbox[t]{4in}{
     \ignorespaces #2}\par }}}
\def\@addressname{}
\def\adr#1{{\normalsize\unskip$^{#1}$}}
\def\@maketitle{%
\def\and{{\rm and}}
  \newpage
  \null
  {\centering
  \let \footnote \thanks
    {\Large\bf   \@title \par}%
    \vskip 1.5em%
      \lineskip .5em%
    {\bf\normalsize   \@author\par}
      \vspace{1em}
    {\small \@addressname}

  }%
  \par
  \vskip 1.5em}
\def\section{\@startsection {section}{1}{\z@}{-3.5ex plus-1ex minus
    -.2ex}{1.5ex plus.2ex}{\reset@font\large\bf}}
\def\subsection{\@startsection{subsection}{2}{\z@}{-3.25ex plus-1ex
    minus-.2ex}{1.5ex plus.2ex}{\reset@font\normalsize\bf}}
\def\subsubsection{\@startsection
     {paragraph}{4}{\z@}{3.25ex plus1ex minus.2ex}{-1em}{\reset@font
     \normalsize\bf}}

\def\cite{\@ifnextchar[{\@tempswatrue\@citex}{\@tempswafalse\@citex[]}}

\def\@citex[#1]#2{%
\if@filesw\immediate\write\@auxout{\string\citation{#2}}\fi
\leavevmode\unskip\ \@cite{\@collapse{#2}}{#1}}

\def\@bylinecite{%
\@ifnextchar[{\@tempswatrue\@CITEX}{\@tempswafalse\@CITEX[]}%
}

\def\@CITEX[#1]#2{%
\if@filesw\immediate\write\@auxout{\string\citation{#2}}\fi
\leavevmode\unskip$^{\scriptstyle\@CITE{\@collapse{#2}}{#1}}$}

\def\@cite#1#2{[{#1\if@tempswa , #2\fi}]} %
\def\@CITE#1#2{{#1\if@tempswa , #2\fi}} %

\def\@collapse#1{%
{%
\let\@temp\relax
\@tempcntb\@MM
\def\@citea{}%
\@for \@citeb:=#1\do{%
\@ifundefined{b@\@citeb}%
{\@temp\@citea{\bf ?}%
\@tempcntb\@MM\let\@temp\relax
\@warning{Citation `\@citeb ' on page \thepage\space undefined}%
}%
{\@tempcnta\@tempcntb \advance\@tempcnta\@ne
\edef\MyTemp{\csname b@\@citeb\endcsname}%
\def\@tempa{\@temptokena=\bgroup}%
\afterassignment\@tempa %
\@tempcntb=0\MyTemp\relax}%
\ifnum\@tempcntb=0\relax%
\@tempcntb=\@MM
\@citea\MyTemp
\let\@temp = \relax
\else %
\edef\@tempd{\number\@tempcntb}%
\ifnum\@tempcnta=\@tempcntb %
\ifx\@temp\relax %
\edef\@temp{\@citea\@tempd}%
\else
\edef\@temp{\hbox{--}\@tempd}%
\fi
\else %
\@temp\@citea\@tempd
\let\@temp\relax
\fi
\fi
}%
\def\@citea{,}%
}%
\@temp %
}%
}%

\catcode `\@=12

\title{  GEOMETRY VIA COHERENT STATES}

\author{Stefan Berceanu}
\address{Institute for Physics and Nuclear Engineering, Department of
Theoretical Physics,
PO BOX MG-6, Bucharest-Magurele, Romania,
E-mail: Berceanu@theor1.ifa.ro}
\newtheorem{com}{Comment}
\newtheorem{thm}{Theorem}
\newtheorem{pr}{Proposition}
\font\ac=eufm10 scaled\magstep1
\font\ab=msbm10 scaled\magstep1
\font\aba=msbm7 
\font\ai=eusb10 scaled \magstep 1
\def\openone{\leavevmode\hbox{\small1\kern-3.8pt\normalsize1}}%
\newcommand{\oo}{\openone}

\newcommand{\gz}{\mbox{\ab Z}}
\newcommand{\gr}{\mbox{\ab R}}
\newcommand{\gcp}{\mbox{\ab CP}}
\newcommand{\gc}{\mbox{\ab C}}
\newcommand{\gcm}{\mbox{\aba C}}
\newcommand{\Gras}{\mbox{$G_n({\gc}^{m+n})$}}
\newcommand{\men}{\mbox{$\widetilde {\bf M}$}} 
\newcommand{\got}[1]{{\mbox{\ac{#1}}}}
\newcommand{\mb}[1]{{\mbox{\boldmath{$#1$}}}}
\newcommand{\gh}{\mbox{\ai H}}
\newcommand{\gl}{\mbox{\ai L}}
\newcommand{\gpl}{\mbox{{\bf P}(\gl )}}
\newcommand{\gph}{\mbox{{\bf P}(\gh )}}

\newcommand{\df}{\mbox{$:=$}}

\begin{document}
\maketitle
\begin{abstract}
It is shown how the coherent states permit to find different geometrical
objects as the geodesics, the conjugate locus, the cut locus,
the Calabi's diastasis and its domain of definition,
the Euler-Poincar\'e characteristic, the number of
Borel-Morse cells,  the Kodaira embedding theorem.
\end{abstract}
\section{Coherent state manifold and coherent vector manifold}
In this talk the coherent states \cite{klauder} are presented as a
 powerful tool
in global differential geometry and algebraic geometry \cite{ber1},\cite{ber2}.
For general references  and proofs see also
\cite{ber6},\cite{ber7}.
The results are illustrated  on the complex Grassmann manifold
 \cite{ber4}.

We start with some notation.
Let  $\pi$ be a unitary irreducible representation,  $G$ a Lie group,  ${\gh}$
a separable complex Hilbert space. Let also the orbit
  $\widetilde {\bf M} ={\widetilde
\pi (G)} \widetilde{\psi}_0~$, where $  \psi_0 \in \gh,
~ \xi :{\gh}\rightarrow {\gph}$ is the projection
$\xi  (\psi_0)= \widetilde{\psi}_0~$. Then there is a diffeomemorphism
$\widetilde {\bf M} \approx G/K $,  where $K$ is the
 stationary group of  $ \widetilde\psi_0$. If
   $\iota :\widetilde {\bf M} \hookrightarrow {\gpl}$ is a biholomorphic
   embedding in some projective Hilbert space,
 then $\widetilde {\bf M}$ is called  {\it coherent state
manifold}.
 If $\sigma:\widetilde {\bf M}\rightarrow{\cal{S}(\gh)}$
  is a local section  in the unit sphere in $\gh$,
 let ${\bf M}' =\sigma(\widetilde {\bf M})$ be the  holomorphic line
bundle associated  to the principal holomorphic
 bundle $P\rightarrow G^c\rightarrow G^c/P$
 by a holomorphic character  $\chi$ of the
parabolic subgroup $ P $  of the compexification  $G^{\gcm}$ of $G$.
 ${\bf M}'$ is a quantization bundle
over \men~ \cite{cgr}.

  The following assertions are equivalent:
there exists the embedding $\iota$; there exists a positive line bundle
 ${\bf M}'$
over $\widetilde{\bf M}$; the line bundle ${\bf M}'$ is ample;
there exists $ m_0$ such  that for $m\geq m_0,~
 {\bf M}={\bf M}'^m=\iota^*[1]$, where $[1]$ is the hyperplane bundle
  over \gpl .

{\bf M} is called {\it coherent vector manifold}.
The  Perelomov's coherent vectors  are
\begin{equation}
{\mb e}_{Z,j}=\exp\sum_{{\varphi}\in\Delta^+_n}(Z_{\varphi}F^+_{\varphi})
j ,~~~~\underline{{\mb e}}_Z=({\mb e}_Z , {\mb e}_Z)^{-1/2}{\mb e}_Z ,
\label{z}
\end{equation}
\begin{equation}
{\mb e}_{B,j}=\exp\sum_{{\varphi}\in\Delta^+_n}(B_{\varphi}F^+_{\varphi}-{\bar
B}_{\varphi}F^-_{\varphi}) j , \label{b}
~~~{\mb e}_{B,j}\df  \underline{{\mb e}}_{Z,j},
\end{equation}
where $j$ is an extreme weight vector,
 $\Delta^+_n$ denotes the positive non-compact roots,
  $ Z\df (Z_ \varphi )
\in {\gc}^d$  are local
coordinates in the neighbourhood ${\cal V}_0 \subset \widetilde {\bf M} $ and
 $ d ={\mbox{\rm dim}}_{\gcm} \widetilde {\bf M}$.
In eqs. (\ref{z}), (\ref{b}) $F^+_{\varphi} j\neq  0,F^-_{\varphi} j = 0,
~  \varphi\in\Delta^+_n $, and $(\underline{{\mb e}}_{ Z^\prime},
\underline{{\mb e}}_{ Z})$ is the    hermitian  scalar
product of holomorphic sections in the holomorphic line bundle ${\bf M}$ over
\men .

\section{Geodesics}

\begin{pr}
 {\it For hermitian symmetric spaces,
  the dependence} $Z(t)=Z(tB)$ {\it appearing when
one passes from} eq. (\ref{b}) {\it to} eq. (\ref{z}) {\it describes in
${\cal V}_0$ a geodesic.}
\end{pr}

For example, on the   Grassmannian  $\Gras =SU(n+m)/S(U(n)\times U(m))$,
 \begin{eqnarray}
\nonumber \Gras & = & \exp\left(\matrix{0&B\cr
                             B^*&0\cr}\right)o=
   \left(\matrix{{\cos}\sqrt{BB^*}&B{\displaystyle
{ {\sin}\sqrt{B^*B}\over \sqrt{B^*B}}}\cr
                                 -{\displaystyle {{\sin} \sqrt{B^*B}\over
                         \sqrt{B^*B}}}B^*&{\cos}\sqrt{B^*B}\cr}\right)o \\
\nonumber &  & \\
\nonumber & = & \left(\matrix{{\oo}&Z\cr 0&{\oo}\cr}\right)
\left(\matrix{({\oo}+ ZZ^*)^{1/2}&0\cr 0&({\oo}+ Z^*Z)^{1/2}\cr}\right)
\left(\matrix{{\oo}&0\cr - Z^*&{\oo}\cr}\right)o~,
\end{eqnarray}
the geodesics in  ${\cal V}_0 $ are  given by the expression
\begin{equation}
Z=B{{\rm tg} \sqrt{B^*B}\over \sqrt{B^*B}} .\label{geo}
\end{equation}

\section{Conjugate locus}

\begin{pr}
 {\it For }
{\it  Hermitian symmetric spaces  $\men$, the  parameters}
 $B_{\varphi}$
 {\it in formula} (\ref{b}) {\it of normalised coherent states are normal 
coordinates
in the normal neighbourhood} ${\cal V}_0$.

\end{pr}
\begin{thm}
 {\it Let  $\widetilde {\bf M}$
 be a Hermitian symmetric space
parametrized in} ${\cal V}_0$  {\it as in eqs.} (\ref{z}),
(\ref{b}). {\it Then the
conjugate locus of the point} $ o $ {\it is obtained vanishing the Jacobian
of the exponential map $Z=Z(B)$ and the corresponding transformations of the
 chart from ${\cal V}_0$}.
\end{thm}

\begin{thm}
\label{sakth}
The tangent conjugate locus $C_0$ of  ${\bf O}\in\Gras$ is
\begin{equation}
\label{ura}
C_0=\bigcup_{k,p,q,i}ad\,k(t_iH)~,~i=1,2,3;~1\leq p<q\leq r,
~ k\in K;
\end{equation}
   \begin{equation}
H=\sum_{i=1}^r h_iD_{i\,n+i},~h_i\in\gr,~\sum h^2_i=1~,
\label{hh}
\end{equation}
\begin{equation}
\begin{array}{l@{\:=\:}c@{\:,\:}l}
t_1  & \displaystyle{\frac{\lambda \pi}{|h_p\pm h_q|}}  &
 ~\mbox{\rm multiplicity}~2;\\[2.ex]
t_2  & \displaystyle{\frac{\lambda \pi}{2|h_p|}} & 
  ~\mbox{\rm multiplicity}~1;\\[2.ex]
t_3  & \displaystyle{\frac{\lambda \pi}{|h_p|}}  & 
 ~\mbox{\rm multiplicity}~2|m-n|; ~\lambda\in \gz^{\star}~.
\end{array}
\label{ttt}
\end{equation}

 The conjugate locus of {\bf O} in \Gras~ is given by the union
\begin{equation}
\label{reun}
{\bf C}_0={\bf C}^W_0\cup {\bf C}^I_0,
\end{equation}

\begin{equation}
\label{ect1}
{\bf C}^I_0= \exp \bigcup_{k,p,q}Ad\,k(t_1H)~,
\end{equation}
\begin{equation}
\label{ect2}
{\bf C}^W_0= \exp \bigcup_{k,p}Ad\,k(t_2H)~.
\end{equation}
 Exponentiating the vectors of the type $t_1H$ we get the points
of ${\bf C}^I_0$ (at least two of the stationary angles with {\bf O}
are equal); $t_2H$ are sent to the points of
 ${\bf C}^W_0$  (at least one of the stationary angles with {\bf O}
is $0$ or $\pi /2$).

 ${\bf C}^W_0$  is given by the disjoint union
\begin{equation}
{\bf C}^W_0=\label{72}
\cases {V^m_1\cup V^n_1,&$ n\leq m,$\cr
            V^m_1\cup V^n_{n-m+1},& $n>m,$ \cr}
\end{equation}
 where
\begin{equation}
V^m_1=\label{74.1}
\cases {{\gcp}^{m-1},& {\mbox{\rm {for}}} $n=1 ,$\cr
              W^m_1\cup W^m_2\cup \ldots W^m_{r-1}\cup W^m_r,&
         {\mbox{\rm {for}}}      $1<n ,$\cr}
\end{equation}
\begin{equation}
W^m_r=\cases {G_r({\gc}^{\max (m,n)}),&$n\not= m,$\cr
             {\bf O}^{\perp},&$n=m ,$\cr}
\end{equation}
\begin{equation}
V^n_1=\label{72.2}
\cases {W^n_1\cup \ldots \cup W^n_{r-1}\cup {\bf O},&$1<n\leq m,$\cr
             {\bf O}, &$n=1 ,$ \cr}
\end{equation}
\begin{equation}
V^n_{n-m+1}=W^n_{n-m+1}\cup W^n_{n-m+2}\cup\ldots\cup W^n_{n-1}\cup
{\bf O}~,~ n>m ~.\label{73.3}
\end{equation}
\end{thm}
\noindent Here $D_{ij}=E_{ij}-E_{ji},$ ${\bf O}^{\perp}$ is the orthogonal
complement of the $n$-plane ${\bf O}$ in ${\gc}^N$ and

\begin{equation}
V^p_l=\left \{Z\in  G_n({\gc}^{n+m})\vert \dim (Z\cap {\gc}^p)
\geq l\right\} ,
\end{equation}
\begin{equation}
W^p_l=V^p_l-V^p_{l+1}=
\left \{Z\in  G_n({\gc}^{n+m})\vert \dim (Z\cap {\gc}^p)
= l\right\} .
\end{equation}

\section{Cut locus}

We remember that $q$ is in the {\it cut locus} ${\bf CL_p}$ of $p\in \men$~
if $q$ is the nearest point to $p$ on the
geodesic emanating from $p$ beyond which the geodesic ceases to minimize his
arc length \cite{kn}.    By {\it polar divisor} of $e_0\in {\bf M}$ we mean
  the set
$\Sigma _0=\left\{ \psi \in {\bf M}|(e_0, \psi)=0\right\}. $

\begin{thm}
 {\it Let} $\widetilde {\bf M}$ {\it be a
 homogeneous manifold} $\widetilde {\bf M}\approx G/K$
{\it
parametrized in the neighbourhood } $ {\cal V}_0$ {\it around $Z=0$ as
in eq.} (\ref{z}). {\it Then } {\it  (the disjoint union)}

\begin{equation}
\widetilde {\bf M} ={\cal V}_0\cup \Sigma _0.\label{reu}
\end{equation}

{\it Moreover, if the condition} $B)$ {\it is true, then}

\begin{equation}
\Sigma _0={\bf CL}_0. \label{clo}
\end{equation}
\end{thm}

\parbox{1cm}{\it A)}\parbox[t]{132mm}{$ {\rm Exp}\vert _o
=\lambda \circ \exp \vert _{\got m}~.  $}

\parbox{1cm}{\it B)}\parbox[t]{132mm}{ On the Lie algebra  {\got g} of $ G$
there exists an $ Ad(G)$-invariant,
symmetric, non-degenerate bilinear form $ B$ such that the restriction
of $B$ to the Lie algebra {\got k}  of $ K$  is likewise non-degenerate.}

\noindent  Here ${\got g}={\got k}\oplus {\got m}$ is the orthogonal
 decomposition of {\got g} with respect to the $B$-form,
${\rm Exp}_p:\widetilde {\bf M}_p\rightarrow \widetilde {\bf M}$ is
 the geodesic exponential map,
 $\exp :{\got g} \rightarrow G $, $o=\lambda (e)$, $e $ is the unit element
  in $G$ and $\lambda$  is the   projection $\lambda :G\rightarrow G/K$.
Note that {\it the symmetric spaces have property $A)$} and if $\men \approx
G/K$ {\it verifies} $B)$, {\it then it also verifies} $A)$ (cf. \cite{kn}).

\begin{com} The cut locus is present everywhere one speaks about
 coherent states.\end{com}

\section{ Calabi's diastasis}

We remember that
the {\it Calabi's diastasis}
is expressed through the coherent states as
 $D(Z',Z)=-2\log\, \vert (\underline{e}_{Z'} ,\underline{e}_{Z})\vert $
 (cf. \cite{cgr}).
 Let also $d_c([\omega '],[\omega ])=\arccos
  \frac{\vert(\omega ',\omega)\vert}{\Vert \omega '\Vert \Vert \omega \Vert }$
denote the hermitian elliptic Cayley distance on the complex projective
space.

\begin{pr}
 The diastasis distance $D(Z',Z)$ between $Z', Z \in {\cal V}_0
 \subset \widetilde {\bf M}$ is related to the geodesic distance
 $\theta = d_c(\iota _*(Z'),\iota _*(Z))$,
 where $\iota :\widetilde {\bf M} \hookrightarrow {\gpl},$ by the relation
$$D(Z',Z)= -2\log\, \cos\,\theta .$$
 
 If $\widetilde {\bf M}_n$ is noncompact, 
  $\iota ':\widetilde {\bf M}_n \hookrightarrow \gcp^{N-1,1}=
  SU(N,1)/S(U(N)\times U(1))$,
and $\delta _n (\theta
  _n)$ is the length of the geodesic joining $\iota '(Z'),\iota '(Z)~$
 (resp. $\iota (Z'),\iota (Z))$,
then $$\cos\,\theta _n=
  (\cosh\,\delta _n)^{-1}=e^{-D/2}.$$
\end{pr} 

\begin{com}
 The relation  (\ref{clo}) furnishes for manifolds of symmetric type
   a geometric description of the
 domain of definition of  Calabi's diastasis: for $z$ fixed, $z'\notin
 {\bf CL}_z$.\end{com}

\section{ The Euler-Poincar\'e characteristic, Borel-Morse cells,
  Kodaira embedding}
 
  \begin{thm}\label{bigtm}
For flag manifolds
 $\widetilde {\bf M} \approx  G/K $, the following 7 numbers are equal:
 
 1) the maximal number of orthogonal coherent vectors;
 
 2) the number of holomorphic global section of the holomorphic line bundle
 $\bf M$;
 
 3) the dimension of the fundamental  representation in the Borel-Weil theorem;

 4) the minimal $N$ appearing in the Kodaira embedding theorem,
 $\iota :\widetilde {\bf M} \hookrightarrow  \gcp^{N-1}$;
 
 5) the number of critical points of the energy function $f_H$ attached to a 
 Hamiltonian $H$ linear in the generators of the Cartan algebra of G,
 with unequal coefficients;

 6) the Euler-Poincar\'e characteristic
$\chi (\widetilde {\bf M})=
 [W_G]/[W_H]$,  $[W_G]= card\,W_G$,  where $ W_G$ denotes the Weyl group of
  $G$;

 7) the number of Borel-Morse cells which appear in the CW-complex
 decomposition of $\widetilde {\bf M}$.
\end{thm}

 \begin{com}
 The Weil prequantization condition is nothing else that
  the condition to have a Kodaira embedding,
 i.e. the algebraic manifold \men~ to be a Hodge one.
\end{com}

\section*{Acknowledgments}
The author expresses his thanks to Professor Heinz-Dietrich
Doebner  for the possibility to attend the  Colloquium.
Discussions during the Colloquium with   Professors
Martin Bordermann, Michael Forger, Joachim Hilgert, John Klauder,
Askold Perelomov and Martin Schlichenmaier are acknowledged.

\end{document}